\documentclass[11pt,twoside]{proc-l}
\usepackage{latexsym}
\usepackage{amssymb,amsbsy,amsmath,amsfonts,amssymb,amscd}

\newtheorem{theorem}{Theorem}[section]

\theoremstyle{definition}
\newtheorem{definition}{Definition}[section]

\theoremstyle{remark}

\numberwithin{equation}{section}

\newcommand{\D}{{\mathcal D}}

\newcommand{\sg}{\sigma}
\newcommand{\Sg}{\beta}
\newcommand{\Chi}{\chi}
\newcommand{\p}{\partial}
\newcommand{\R}{{\mathbb R}}

\newcommand{\TP}{{\mathbb T_P}}
\newcommand{\fer}[1]{(\ref{#1})}

\newcommand{\be}{\begin{equation}}
\newcommand{\ee}{\end{equation}}

\def\<#1>{\left\langle #1\right\rangle}

\begin{document}

\title[Anisotropic Degenerate Parabolic-Hyperbolic Equations]{Large-Time Behavior
of Periodic Entropy Solutions to Anisotropic Degenerate
Parabolic-Hyperbolic Equations}

\author{Gui-Qiang Chen}
\address{G.-Q. Chen: School of Mathematical Sciences, Fudan University,
 Shanghai 200433, China; Department of Mathematics, Northwestern
University, Evanston, IL 60208-2730, USA.}
\email{gqchen@math.northwestern.edu}
\thanks{Gui-Qiang Chen's research was supported in
part by the National Science Foundation under Grants DMS-0807551,
DMS-0720925, and DMS-0505473, and the Natural Science Foundation of
China under Grant NSFC-10728101. This paper was written as part of
the international research program on Nonlinear Partial Differential
Equations at the Centre for Advanced Study at the Norwegian Academy
of Science and Letters in Oslo during the academic year 2008--09.}

\author{Beno\^ \i t Perthame}
\address{B. Perthame: UPMC, Univ. Paris 06, UMR 7598, Laboratoire Jacques-Louis Lions,
4, pl. Jussieu, 75005, Paris, France, and Institut Universitaire de
France}
\email{benoit.perthame@upmc.fr}
\thanks{}

\date{September 8, 2008}
%
\begin{abstract}
We are interested in the large-time behavior of periodic entropy
solutions in $L^\infty$ to anisotropic degenerate
parabolic-hyperbolic equations of second-order. Unlike the pure
hyperbolic case, the nonlinear equation is no longer self-similar
invariant and the diffusion term in the equation significantly
affects the large-time behavior of solutions; thus the approach
developed earlier based on the self-similar scaling does not
directly apply. In this paper, we develop another approach for
establishing the decay of periodic solutions for anisotropic
degenerate parabolic-hyperbolic equations. The proof is based on
the kinetic formulation of entropy solutions. It involves time
translations and a monotonicity-in-time property of entropy solutions, and
employs the advantages of the precise kinetic equation for the
solutions in order to recognize the role of
nonlinearity-diffusivity of the equation.
\end{abstract}
\maketitle

\section{Introduction}

We study the large-time behavior of periodic solutions in $L^\infty$
to nonlinear {\em anisotropic} degenerate parabolic-hyperbolic
equations of second-order. Consider the Cauchy problem for the second-order equations:
\begin{equation}
\p_t u+\nabla_x\cdot f(u)=\nabla_x\cdot (A(u)\nabla_x u), \qquad
x\in \R^d, \quad t\geq 0, \label{PSCL1}
\end{equation}
\begin{equation}
u|_{t=0}=u_0 \in L^\infty(\R^d), \qquad u_0(x+P)=u_0(x) \quad
\text{a.e.} \; x\in \R^d, \label{PSCL2}
\end{equation}
where  $P=(P_1,...,P_d)$ is the period with $ P_i>0$,
$\TP:=\Pi_{i=1}^d[0, P_i]$, $f: \R\to \R^d$ satisfies
\begin{equation}
a(\cdot): =f^\prime (\cdot) \in L^\infty_{\rm loc}(\R; \R^d),
\label{Ass1}
\end{equation}
and the $d\times d$ matrix $A(u)=(a_{ij}(u))$ is symmetric,
nonnegative, and locally bounded so that it can be always written
under the form
\begin{equation}
a_{ij}(u)= \sum\limits_{k=1}^{d} \sg_{ik}(u) \sg_{jk}(u), \qquad
\sg_{ik} \in L^\infty_{\rm loc}(\R), \label{Ass2}
\end{equation}
with $(\sigma_{ik}(u))$ the square root matrix of $A(u)$.

\smallskip
Equation \eqref{PSCL1} and its variants model degenerate
diffusion-convection motions of ideal fluids and arise in a wide
variety of important applications (cf. \cite{BLS,BCBT,CJ,EFM,npp}
and the references cited therein), for which a deep understanding of
solutions to \eqref{PSCL1} is in great demand, at both the
theoretical and numerical level.

\smallskip
In \cite{ChenPerthame}, a well-posedness theory has
been established for $L^{1}$ solutions of the Cauchy problem
\eqref{PSCL1} and \eqref{PSCL2} of anisotropic degenerate
parabolic-hyperbolic equations of second-order.
It extends the isotropic theory for degenerate parabolic-hyperbolic
equations and  several latter studies, see for instance
\cite{Ca,ChenDiBenedetto,KarlsenR,MV} and the references therein.
A notion of kinetic solutions, a new concept in this context, and a
corresponding kinetic formulation have been extended. In particular,
it has been also proved that, when $u_0\in L^\infty$, the kinetic
solution is equivalent to the entropy solution, or to the
dissipative solution \cite{PeSou}, which is unique; this provides a
new path to study the behavior of entropy solutions in $L^\infty$
through the corresponding kinetic equations for anisotropic
degenerate parabolic-hyperbolic equations. In this paper, we employ
the advantages of this path to develop a new approach for
establishing the decay of periodic entropy solutions in $L^\infty$
to \eqref{PSCL1} and \eqref{PSCL2} as $t\to\infty$. The main theorem
of this paper is the following.

\begin{theorem}[Main theorem]\label{thm:decay}  Assume that {\rm \fer{Ass1}}
and {\rm \fer{Ass2}} hold. Let $u\in L^\infty
([0,\infty)\times\R^d)$ be the unique periodic entropy solution to
{\rm \fer{PSCL1}} and {\rm \fer{PSCL2}}. Assume that the flux
function and the diffusion matrix $A(u)$ satisfy the
nonlinearity-diffusivity condition: For any $\delta>0$,
\begin{equation}\label{nonlinearity}
\sup_{|\tau|+ |\kappa| \geq \delta} \; \int_{|\xi| \leq
\|u^0\|_\infty} \frac{ \lambda \; d\xi}{\lambda +|\tau +a(\xi)\cdot
\kappa|^2+ (\kappa^\top A(\xi)\kappa)^2} := \omega_\delta(\lambda)
{\;}_{\overrightarrow{\; \lambda \rightarrow 0 \; }}\;  0 .
\end{equation}
 Then we have
\begin{equation}\label{decay}
\|u(t,\cdot)-\bar{u} \|_{L^1(\TP)} {\;}_{\overrightarrow{\; t \rightarrow \infty \; }}\;  0,
\end{equation}
where
$$
\bar{u}=\frac{1}{|\TP|}\int_\TP u_0(x)dx.
$$
\end{theorem}

The {\em nonlinearity-diffusivity condition} \eqref{nonlinearity}
for equation {\rm \fer{PSCL1}} is developed from \cite{GLPS,LPT2}
and is reminiscent from the theory of velocity averaging lemmas
\cite{GLPS,Pe2,PeSou1,TT}.
For smooth coefficients, condition (\ref{nonlinearity}) is
equivalent to the simpler and more standard setting: For any
$(\tau,\kappa)\in \R^{d+1}$ with $\tau^2+|\kappa|^2=1$, we have
$$
\mathcal{L}^1\big\{\xi \in \R\,:\, |\xi|\leq \|u_0\|_\infty, \; \tau
+a(\xi)\cdot \kappa=0, \,\, \displaystyle \kappa^\top
A(\xi)\kappa=0\big\}=0.
$$
Several explicit examples are given in \cite{LPT2,TT}. It implies
that there is no interval of $\xi$ on which both the flux function
$f(\xi)$ is affine and the diffusion matrix $A(\xi)$ is degenerate,
and thus also makes the relation with the applications of the theory
of compensated compactness \cite{Mu,Ta} to one-dimensional
hyperbolic conservation laws. Such a nonlinearity-diffusivity
condition is necessary for the decay of periodic solutions and the
compactness of solution operators. See also the books
\cite{dafermos,serre}.

Theorem 1.1 also extends naturally the non-degeneracy condition for
the purely hyperbolic case in \cite{EE} where the first long-time
convergence result of periodic solutions was obtained in one or two
dimensions for $BV$ initial data and higher order local
non-degeneracy conditions that replace (\ref{nonlinearity}).

\smallskip
Unlike the pure hyperbolic case, equation {\rm \fer{PSCL1}} is no
longer self-similar invariant and the diffusion term in the equation
significantly affects the large-time behavior of solutions; thus the
approach in \cite{ChenFrid} based on the self-similar scaling for
the pure hyperbolic case does not apply and we have to change the
strategy of proof. The approach developed in this paper is based on
the kinetic formulation of entropy solutions in \cite{ChenPerthame},
involves time translations and a monotonicity-in-time of entropy
solutions, and employs the advantages of the kinetic equations of
the solutions, in order to recognize the role of the
nonlinearity-diffusivity of the equation.

\smallskip
The rest of this paper is organized as follows. We first recall the
notion of entropy solutions and their precise kinetic formulation,
and then analyze some basic properties of entropy solutions.
Finally, we develop a new approach to give a rigorous proof for the
long-time asymptotic result.

\section{Entropy Solutions and Kinetic Formulation}

In this section, we first recall the notion of entropy solutions and
their precise kinetic formulation, which requires some care to
define appropriately the various terms of the equation. Then we
analyze some basic properties of entropy solutions, which will be
used in the proof of Theorem \ref{thm:decay}.

\begin{definition} \label{Def:1}
An entropy solution is a function
$u(t,x)\in L^\infty\big([0,\infty)\times \R^d\big)$ such that

\noindent {\rm (i)}
$\sum\limits_{i=1}^{d} \p_{x_i}\Sg_{ik}(u)\in L^2 \big(
[0,\infty)\times \R^d\big)$,
$k=1,\cdots,d$, for $\beta_{ik}(u)=\int^u\sigma_{ik}(v)\, dv$;

\noindent {\rm (ii)} For any function $\psi\in C_0(\R)$ with
$\psi(u)\geq 0$ and any $k=1,\cdots,d$, the chain rule holds:
\begin{equation}
\sum\limits_{i=1}^{d}\p_{x_i}\Sg^\psi_{ik}(u)
=\sqrt{\psi(u)}\sum\limits_{i=1}^{d}\p_{x_i}\Sg_{ik}(u) \in L^2
([0,\infty)\times \R^d), \; \text{ for } \;
(\Sg^{\psi}_{ik})^\prime= \sqrt{\psi}\Sg_{ik}^\prime; \label{DES0}
\end{equation}

\noindent {\rm (iii)} For any smooth function $S(u)$, there exists
an entropy dissipation measure $m^{S''}(t,x)$ satisfying that
\begin{equation}\label{DES3}
m^{S''}(t,x) =\int_{\R} S''(\xi) \, m(t,x,\xi) \, d\xi \quad
\hbox{with}\,\, m(t,x,\xi) \, \hbox{a nonnegative measure,}
\end{equation}
such that
\begin{equation}
\p_t S(u) +\nabla_x\cdot
q^S(u)
-\nabla_x\cdot\big(A(\xi)\nabla_xS(u)\big)
= -(m^{S''} + n^{S''}) \label{DES2}
\end{equation}
in ${\D}^\prime (\R^+\times \R^d)$ with initial data
$S(u)|_{t=0}=S(u_0)$, where $q^S: \R\to\R^d$ is the corresponding
entropy flux, i.e. $(q^S(u))'=S'(u)f'(u)$, and
$$
n^{S''}(t,x,\xi)=\int S''(\xi)n(t,x,\xi)d\xi
$$
for $n(t,x,\xi)$ the
parabolic defect measure of $u(t,x)$ defined as:
\begin{equation}
n(t,x,\xi) :=\delta(\xi-u(t,x)) \sum\limits_{k=1}^{d} \big(
\sum\limits_{i=1}^{d}
    \p_{x_i}\Sg_{ik}(u(t,x))\big)^2  \quad a.e.
\label{DES1}
\end{equation}
\end{definition}

We point out that the chain rule in (ii) has to be assumed in the
anisotropic case, and this makes the main difference with the
isotropic case in \cite{Ca} where this property follows from an
argument reminiscent to the theory of Sobolev spaces. The
requirement (iii) has been made with notations that are adapted to
the kinetic formulation we introduce now. The $L^2$-condition in (i)
is required to define the parabolic defect measure $n(t,x,\xi)$ in
\eqref{DES1} (also see \cite{ChenDiBenedetto}).

To do so, we may factor out an $S'(u)$ in equation \fer{DES2} and
obtain a more handful kinetic formulation of nonlinear degenerate
parabolic-hyperbolic equations of second-order with form
\fer{PSCL1}. The new ingredient of this formulation is the
identification of the kinetic defect measure $m(t,x,\xi)$ and the
degenerate parabolic defect measure $n(t,x,\xi)$ in a precise
manner, even in the region where $u(t,x)$ is discontinuous. Compare
with the classical {\it kinetic formulation} for scalar hyperbolic
conservation laws in \cite{LPT2} (see also \cite{Pe2}).

\medskip
We introduce the kinetic function $\Chi$ on $\R^2$:
\begin{equation}
     \Chi(\xi;u)=\left \{
      \begin{array}{ll}
      +1 &\ \ {\rm for}\ \ 0 < \xi < u, \\
      -1 &\ \ {\rm for}\ \ u < \xi < 0,\\
      \;0 &\ \ {\rm otherwise}.
      \end{array}
      \right.
\label{PDefChi}
\end{equation}
We notice that, if $u$ is in $L^\infty\big([0,\infty)\times
\R^d\big)$ and periodic in $x$ with period $P$, then $\Chi(\xi;u)$
is in $L^\infty\big([0,\infty)\times\R^{d+1}\big)$ and periodic in
$x$ with period $P$.

\medskip
The simple representation formula $ S(u)=\int_\R S'(\xi)\, \Chi(\xi;u) \,
d\xi $ leads to the following kinetic equation, which is equivalent
to the entropy identity  \fer{DES2}:
\begin{equation}
\p_t\Chi(\xi;u) +a(\xi)\cdot \nabla_x \Chi(\xi;u)
-\nabla_x\cdot\big(A(\xi)\nabla_x \chi(\xi;u)\big)
=\p_\xi (m+ n)(t,x,\xi) \label{KS1}
\end{equation}
in ${\D}^\prime (\R^+\times \R^{d+1})$ with initial data
\begin{equation}\label{KS2}
\Chi(\xi;u)|_{t=0} = \Chi(\xi;u_0),
\end{equation}
where $n(t,x,\xi)$ is still defined through \fer{DES1}.

In \cite{ChenPerthame}, it has been proved that the entropy
solutions in $L^\infty$ are equivalent to the {\it kinetic
solutions} determined by the kinetic formulation
\eqref{PDefChi}--\eqref{KS2}.
Furthermore, we have

\begin{theorem}\label{thm:periodic}  Assume that {\rm \fer{Ass1}}
and {\rm \fer{Ass2}} hold. Then
\\
{\rm (i)} there exists a unique entropy solution $u\in L^\infty
([0,\infty)\times \R^d)$ to {\rm \fer{PSCL1}}--{\rm \fer{PSCL2}}
such that $u(t, x+P)=u(t,x)$ a.e.
\\
{\rm (ii)} the entropy solution $u(t,x)$ satisfies

{\rm (a)} $\|u(t, \cdot)\|_{L^\infty}\le \|u_0\|_{L^\infty}$;

{\rm (b)} for any $t_2>t_1>0$,
\begin{eqnarray}
&&\int_\TP  |u(t_2,x)|^2\, dx\le \int_\TP |u(t_1,x)|^2\,
dx, \label{monotone1}\\
&&\int_\TP |u(t_2,x)-v|\, dx\le \int_\TP |u(t_1,x)-v|\,
dx\qquad\text{for any constant $v$}. \label{monotone2}
\end{eqnarray}
{\rm (iii)} Moreover, if the flux function $f(u)$ and the diffusion
matrix $A(u)$ further satisfy the nonlinearity-diffusivity condition
\eqref{nonlinearity}, then the solution operator $
u(t,\cdot)=S_tu_0(\cdot)\,:\, L^\infty\to L^1$ is locally compact in
$(t,x)$ for $t>0$.
\end{theorem}

The results {\rm (i)}--{\rm (ii)} in Theorem 2.1 are direct
corollaries of the well-posedness results and the arguments in
 \cite{ChenPerthame} (i.e. standard entropy inequalities for
$|u-v|$ and $u^2$) as in the hyperbolic case \cite{dafermos, serre};
while (iii) is a direct corollary of the
kinetic averaging compactness result of \cite{LPT2}; also see
more recent results in \cite{JP,PeSou1, TT}.

\section{Decay of Periodic Entropy Solutions: Proof of Theorem \ref{thm:decay}}

In this section, we develop a new approach to give a rigorous proof
for the decay property of periodic solutions, which takes the
advantage of the precise kinetic formulation
\eqref{PDefChi}--\eqref{KS2}. Without loss of generality, we first
set $\int_\TP u(t,x) dx=0$; otherwise, we may replace $u(t,x)$ by
$u(t,x)-\bar{u}$, $f(u)$ by $f(u+\bar{u})$, and $A(u)$ by
$A(u+\bar{u})$ in \eqref{PSCL1}, so that all the arguments below
remain unchanged. Then, we divide the proof into four steps.

\medskip {\em 1. Limit}.  Theorem 2.1 indicates that the periodic
solution $u(t,x)$ belongs to $L^\infty$, bounded by
$\|u_0\|_{L^\infty}$, and is compact as the solution operator. Also the
function
$$
I(t):=\int_\TP  |u(t,x)|^2\, dx
$$
is a non-increasing, bounded function, which implies that the
following limit exists:
\begin{equation}\label{3.2a}
\lim_{t\to \infty}I(t)=I(\infty)=:I_\infty \in [0, \infty).
\end{equation}

\smallskip
{\em 2. Translations}. Set
$$
v_k(t,x):=u(t+k, x).
$$
Then we find that, for $t\geq -k$,

(i) $\|v_k(t,\cdot)\|_{L^\infty}=\|u(t+k, \cdot)\|_{L^\infty}
    \le \|u_0\|_{L^\infty}$;

(ii) $v_k(t,x)$ is also a periodic entropy solution to {\rm
\fer{PSCL1}}--{\rm \fer{PSCL2}};

(iii) for each $k>0$, $\Chi(\xi; v^k(t,x))$ satisfies
\begin{equation}\label{kinetic:eq-id}
\p_t\Chi(\xi; v^k) +a(\xi)\cdot \nabla_x \Chi(\xi;v^k)
-\nabla_x\cdot\big(A(\xi)\nabla_x\Chi(\xi;v^k)\big)
=\p_\xi (m^k+ n^k)(t,x,\xi)
\end{equation}
in ${\D}^\prime ((-k, \infty)\times \R^{d+1})$.

\smallskip
By Theorem 2.1 (iii) applied to $v^k$, there exists a subsequence
$\{v^{k_j}\}_{j=1}^\infty \subset\{v^k\}_{k=1}^\infty$ and
$v(t,x)\in L^\infty\big(\R^{d+1}\big)$, with $\int_\TP v(t,x)dx=0$,
such that
$$
v^{k_j}(t,x)\to v(t,x) \qquad \text{a.e.}\,\, (t,x)\in \R^{d+1}
\qquad\mbox{as}\,\, j\to \infty.
$$
Correspondingly, we have
$$
\Chi(\xi; v^k(t,x))\to \Chi(\xi; v(t,x)) \qquad \text{a.e.}\,\,
(t,x, \xi) \qquad\mbox{as}\,\, j\to \infty.
$$

\smallskip
Furthermore, multiplying both sides of \eqref{kinetic:eq-id} by
$\xi$ and then integrating over $(t,x,\xi)\in (-T,T)\times
\TP\times\R$ for any $T>0$, we obtain
\begin{equation}\label{measure-bound}
\int_{-T}^T\int_\TP (m^{k}+n^{k})(t,x,\xi)\, dtdxd\xi \le
\frac{1}{2}\big(I(k-T)-I(k+T)\big)\le
\frac{1}{2}|\TP|\|u_0\|^2_{L^\infty}.
\end{equation}
This implies that the nonnegative measure sequence
$(m^{k}+n^{k})(t,x,\xi)$ is uniformly bounded  in $k$ over $(-T,T)\times
\TP\times\R$, and hence there exists a subsequence $k_j$ and
a measure $M(t,x,\xi)$ such that
$$
(m^{k_j}+n^{k_j})(t,x,\xi)\,\, {\rightharpoonup}\,\, M(t,x,\xi)\ge 0
\qquad\mbox{weakly in}\,\, \mathcal{M} \quad\text{as}\,\,
j\to\infty.
$$

\smallskip
On the other hand, since $I(t)$ converges, we also have
\begin{equation}\label{I-zero}
I(k-T)-I(T+k)\to 0 \qquad \mbox{as}\,\, k \to \infty.
\end{equation}
Then we conclude from \eqref{measure-bound} that
$M((-\infty,\infty)\times \TP\times\R)=0$, which implies
\begin{equation}\label{measure-zero}
M(\R^{d+2})=0.
\end{equation}

\smallskip
Furthermore, letting $k\to \infty$ in \eqref{kinetic:eq-id}, we
conclude that $\Chi(\xi; v)$ is a $\TP$-periodic solution in
${\D}^\prime (\R\times \R^{d+1})$ of
\begin{equation}\label{Chi-eq}
\p_t\Chi +a(\xi)\cdot \nabla_x \Chi
-\nabla_x\cdot\big(A(\xi)\nabla_x\Chi\big)=0.
\end{equation}

In particular, multiplying \eqref{Chi-eq} by $\xi$ and then
integrating $dx d\xi$, we have
\begin{equation}
\int_\TP |v(t,x)|^2 dx = I_\infty\in [0, \infty), \qquad \forall\, t
\in \R, \label{eq:cst}
\end{equation}
where $I_\infty=I(\infty)$ is a constant, independent of $t$,
determined in \eqref{3.2a}.

The rest of the proof consists in showing that such a function
$\Chi$ is very particular and is in fact constant (see also this
type of ``rigidity" results in \cite{JOP, DOW,DaPe}).

\medskip
{\em 3. Claim}: $v(t,x)\equiv 0$ a.e. for $(t,x)\in\R\times\R^d$.
%
%

\smallskip
We introduce a time truncation function $\phi(t), 0\le\phi(t)\le 1$,
so that $\phi \Chi$ belongs to $L^2(\R\times\TP\times\R)$. Then we
have
\begin{equation}\label{Chi-eqt}
\p_t(\phi \Chi) +a(\xi)\cdot \nabla_x (\phi \Chi)
-\nabla_x\cdot\big(A(\xi)\nabla_x (\phi \Chi) \big)=
\Chi \, \p_t\phi \qquad \text{in ${\D}^\prime (\R^{d+2})$}.
\end{equation}

\smallskip
Next, since $\phi\Chi$ and $\Chi\phi_t$ are periodic in $x$, we take
the global Fourier transform in $t\in\R$ and the local Fourier
transform in $x\in\TP$ to obtain $\hat{g}(\tau, \kappa; \xi)$ for
$(\phi \Chi)(t,x,\xi)$ and $\hat{h}(\tau, \kappa; \xi)$ for $(\Chi\,
\p_t\phi)(t,x,\xi)$ in $L^2$, where the frequencies
$\kappa=(\kappa_1, \cdots, \kappa_d)$ are discrete:
$$
\kappa_i\frac{2\pi}{P_i}, \qquad n=0, \pm 1, \pm 2, \cdots.
$$
That is, for example,
$$
\hat{g}(\tau, \kappa; \xi)=\frac{1}{|\TP|} \int_\R\int_\TP (\phi
\Chi)(t,x,\xi)\, e^{-i(\tau t+ \kappa\cdot x)}dtdx
$$
so that
$$
(\phi\Chi)(t,x,\xi)=\frac{1}{2\pi}\sum_{\kappa}
\int_\R\hat{g}(\tau,\kappa;\xi)\, e^{i(\tau t+ \kappa\cdot x)}d\tau.
$$

Taking the global Fourier transform in $t\in\R$ and the local
Fourier transform in $x\in\TP$ in both sides of \eqref{Chi-eqt}, we
obtain
$$
\Big(i\big(\tau+a(\xi)\cdot \kappa\big)+ \kappa^\top A(\xi)\kappa
\Big)\hat{g}=\hat{h}.
$$
Following usual ideas from the kinetic averaging lemmas, we may introduce
a free parameter $\lambda >0$ (to be chosen later on) and write
$$
\Big(\sqrt{\lambda} + i\big(\tau+a(\xi)\cdot \kappa\big)+
\kappa^\top A(\xi)\kappa
\Big)\hat{g}=\hat{h} + \sqrt{\lambda}\hat{g}.
$$
This leads to
$$
\hat{g} =( \hat{h} +  \sqrt{\lambda}\hat{g})\;
\frac{1}{\sqrt{\lambda} + i\big(\tau+a(\xi)\cdot
\kappa\big)+\kappa^\top A(\xi)\kappa}.
$$
Integrating in $\xi$ and using the Cauchy-Schwarz inequality, we
find
\begin{eqnarray*}
&&|\widehat{\phi v}|^2(\tau,\kappa)\\
&&\leq 2 \Big(\int_\R \hat{h}^2 d\xi +\lambda \int_\R \hat{g}^2 d\xi
\Big) \,\, \int_\R \Big|\frac{1}{\sqrt{\lambda}+i
\big(\tau+a(\xi)\cdot\kappa\big) +\kappa^\top A(\xi)\kappa}
\Big|^2  d\xi.
\end{eqnarray*}

Notice that the frequencies $\kappa$ are discrete and
may include $\kappa=0$. In particular, when $\kappa\ne 0$, then
there exists $\delta_0>0$ such that $|\kappa|\ge \delta_0$.
On the other hand, since $v(t,x)$ has mean
zero in $x$ over $\TP$, we have
$$
\widehat{\phi v}(\tau,0)=0.
$$

The nonlinearity-diffusivity condition (\ref{nonlinearity}) gives
that, when $\kappa \neq 0$, for any $\delta\in (0, \delta_0)$,
$$
|\widehat{\phi v}|^2 \leq C \frac {\omega_\delta(\lambda)}{\lambda}
\int_\R \hat{h}^2 d\xi + C \omega_\delta(\lambda)   \int_\R
|\hat{g}|^2 d\xi,
$$
and thus
\begin{eqnarray*}
&&\sum_{\kappa\neq 0} \int_\R| \widehat{\phi v}|^2 d\tau \\
&&\leq C\frac{\omega_\delta(\lambda)}{\lambda} \sum_{\kappa\neq
0}\int_{\R^2} \hat{h}^2 d\xi d\tau + C \omega_\delta(\lambda)
\sum_{\kappa\neq 0}  \int_{\R^2} |\hat{g} |^2
d\xi d\tau\\
&&\leq C \frac{\omega_\delta(\lambda)}{\lambda}
\int_{\R\times\TP\times\R} (\Chi \phi_t)^2 dt d xd\xi + C
\omega_\delta(\lambda)   \int_{\R\times\TP\times\R} |\phi \chi |^2
dt d xd\xi.
\end{eqnarray*}
This implies
$$ \int_{\R\times \TP}| \phi v |^2 dt d x
 \leq C\frac{\omega_\delta(\lambda)}{\lambda}
 \int_{\R\times\TP} |\phi_t|^2 |v| dt dx+
 C\omega_\delta(\lambda)  \int_{\R\times\TP} |\phi|^2 |v| dt dx.
$$
Using (\ref{eq:cst}) and the Cauchy-Schwartz inequality, we arrive
at
\begin{eqnarray}
&\qquad I_\infty \int_{\R}| \phi  |^2 dt
 &\leq C \omega_\delta(\lambda) \Big(\int_{\TP} |v|^2 dx\Big)^{1/2}
 \Big(\frac{1}{\lambda}\int_{\R}|\phi_t|^2 dt + \int_{\R}|\phi|^2 dt\Big)
 \label{3.8-a}\\
&& \leq C  I_\infty^{1/2} \omega_\delta(\lambda)
   \Big(\frac{1}{\lambda}\int_{\R}|\phi_t|^2 dt + \int_{\R} |\phi|^2
   dt\Big).\nonumber
\end{eqnarray}
Choosing first $\lambda$ small and then $\int_\R |\phi'(t)|^2$
small, we conclude from \eqref{3.8-a} that
$$
I_\infty=0,
$$
which implies from (\ref{eq:cst}) that
\begin{equation}\label{3.9-a}
v(t,x)\equiv 0 \qquad \text{a.e.} \,\, (t,x)\in\R\times\R^d.
\end{equation}

\medskip
On the contrary, if $I_\infty>0$, then we can choose $\lambda$ small
enough so that $C \omega_\delta(\lambda)/I_\infty^{1/2} \le
\frac{1}{2}$ and find from \eqref{3.8-a} that
$$
 I_\infty\int_{\R}| \phi  |^2 dt \leq
2CI_\infty^{1/2}
\frac{\omega_\delta(\lambda)}{\lambda}\int_{\R}|\phi_t|^2 dtdx d\xi.
$$
It remains to choose a sequence of functions $\phi_B(t) = 1$ for
$|t|\leq B$, with $B$ a given large number and
$\phi_B'(t)= \frac{2B-|t|}{B}$ for $B\leq |t|\leq 2B$,
and $\phi_B(t)=0$ for $|t|\geq 2B$. In the above
inequality, we find
$$
I_\infty^{1/2} \leq  C \frac{\omega_\delta(\lambda)}{B^2\lambda},
$$
where $C>0$ is a constant independent of $B$ and $\lambda$. When $B$
tends to $\infty$, this implies that $I_\infty$ must vanish, which
is a contradiction.

Therefore, \eqref{3.9-a} holds.

\medskip
{\em 4. Conclusion}. Then
$$
\int_0^1\int_\TP |v^{k_j}(s,x)|\, dxds \to 0 \qquad \text{as}\,\,
j\to \infty.
$$
Therefore, for any $T>k_j+1>k_j+s$ for $s\in (0,1)$, we employ
\eqref{monotone1} for the monotonicity-in-time of solution to obtain
$$
\int_0^1\int_\TP |v^{k_j}(t,x)|\, dxdt \ge
\int_\TP |u(1+k_j,x)|\, dx\ge \int_\TP |u(T,x)|\, dx.
$$

We conclude that
$$
\int_\TP  |u(T,x)|\, dx \to 0\qquad\mbox{as}\,\,\, T\to\infty.
$$
This completes the proof.

\begin{thebibliography}{10}




\bibitem{BLS} M. Bendahmane, M. Langlais, and M. Saad,
\textit{On some anisotropic reaction-diffusion systems with $L\sp
1$-data modeling the propagation of an epidemic disease}, Nonlinear
Anal.  \textbf{54} (2003),  no. 4, 617--636.

\bibitem{BCBT} M.~C. Bustos, F. Concha, R. B\"{u}rger, and E.~M. Tory,
\textit{Sedimentation and Thickening: Phenomenological Foundation
and Mathematical Theory}, Kluwer Academic Publishers: Dordrecht,
Netherlands, 1999.

\bibitem{Ca} J. Carrillo, \textit{Entropy solutions for nonlinear degenerate
problems},
 Arch. Rational Mech. Anal. \textbf{147} (1999), 269--361.

\bibitem{CJ} G. Chavent and J. Jaffre,
\textit{Mathematical Models and Finite Elements for Reservoir
Simulation}, North Holland: Amsterdam, 1986.

\bibitem{ChenDiBenedetto} G.-Q. Chen and E. DiBenedetto,
\textit{Stability of entropy solutions to the Cauchy problem for a
class of hyperbolic-parabolic equations},
     SIAM J. Math. Anal. \textbf{33} (2001), 751--762.

\bibitem{ChenFrid} G.-Q. Chen and H. Frid,
\textit{Decay of entropy solutions of nonlinear conservation laws},
{Arch. Rational Mech. Anal.} \textbf{146(2)} (1999), 95--127.

\bibitem{ChenPerthame} G.-Q. Chen and B. Perthame,
\textit{Well-posedness for nonisotropic degenerate
parabolic-hyperbolic equations}, Annales de l'Institut Henri
Poincar\'{e}: Analyse Non Lin\'{e}aire, \textbf{20} (2003),
645--668.

\bibitem{dafermos} C. M. Dafermos,  \textit{Hyperbolic Conservation Laws
in Continuum Physics}, Second edition,
Springer-Verlag: Berlin, 2005.

\bibitem{DaPe} A.-L. Dalibard and B. Perthame,
\textit{Existence of solutions of the hyperbolic Keller-Segel
model}, Trans. Amer. Math. Soc., to appear.

\bibitem{DOW}  C. De Lellis, F. Otto, and M. Westdickenberg,
\textit{Structure of entropy solutions for multi-dimensional scalar
conservation laws}, Arch. Rational Mech. Anal. \textbf{170} (2003)
137--184.

\bibitem{EE} B.  Engquist and W. E,
\textit{Large time behavior and homogenization of solutions of
two-dimensional conservation laws}, Comm. Pure Appl. Math.
\textbf{46} (1993), 1--26.

\bibitem{EFM} M.~S. Espedal, A. Fasano, and A. Mikeli\'c,
\textit{Filtration in Porous Media and Industrial Applications},
Lecture Notes in Math. \textbf{1734}, Springer-Verlag: Berlin, 2000.

\bibitem{GLPS} F. Golse, P.-L. Lions, B. Perthame, and R. Sentis,
\textit{Regularity of the moments of the solution of a transport
equation}, J. Funct. Anal. \textbf{76} (1988), 110--125.

\bibitem{JOP} P.-E. Jabin, F. Otto, and B.  Perthame,
\textit{Line-energy Ginzburg-Landau models: zero-energy states},
Ann. Sc. Norm. Super. Pisa Cl. Sci. \textbf{(5) 1} (2002), 187--202.

\bibitem{JP} P.-E. Jabin and B. Perthame,
\textit{Regularity in kinetic formulations via averaging lemmas}. A
tribute to J. L. Lions, ESAIM Control Optim. Calc. Var. \textbf{8}
(2002), 761--774 (electronic).

\bibitem{KarlsenR}
K.~H. Karlsen and N.~H. Risebro, \textit{Convergence of finite
difference schemes for viscous and inviscid
  conservation laws with rough coefficients},
M2AN Math. Model. Numer. Anal. \textbf{35(2)} (2001), 239--269.


\bibitem{LPT2} P.-L. Lions, B. Perthame, and E. Tadmor,
\textit{A kinetic formulation of multidimensional scalar
conservation laws and related equations}, {J. Amer. Math. Soc.}
\textbf{7} (1994), 169--191.

\bibitem{MV} A. Michel and J. Vovelle,
\textit{Entropy formulation for parabolic degenerate equations with
general Dirichlet boundary conditions and application to the
convergence of FV methods}, SIAM J. Numer. Anal. \textbf{41} (2003),
2262--2293 (electronic).

\bibitem {Mu} F. Murat, \textit{Compacit\'e par compensation}, Ann. Sc. Norm.
Sup. Pisa, \textbf{5} (1978), 489--507.

\bibitem{npp} J. Nolen,  G. Papanicolaou and O. Pironneau,
\textit{A framework for adaptive multiscale methods for elliptic
problems}, Multiscale Model. Simul. \textbf{7} (2008), 171--196.

\bibitem{Pe2} B. Perthame, \textit{Kinetic Formulations of
Conservation Laws}, Oxford Univ. Press: Oxford, 2002.

\bibitem{PeSou1} B.  Perthame and P.~E. Souganidis,
\textit{A limiting case for velocity averaging}, Ann. Sci. \'Ecole
Norm. Sup. (4) \textbf{31} (1998), 591--598.

\bibitem{PeSou} B. Perthame and P.~E. Souganidis,
\textit{Dissipative and entropy solutions to non-isotropic
degenerate parabolic  balance laws},  Arch. Rational Mech. Anal.
\textbf{170} (2003), 359--370.

\bibitem{serre} D. Serre,  \textit{Systems of Conservation Laws},
Cambridge University Press: Cambridge, 2000.

\bibitem{TT} E. Tadmor and T. Tao, \textit{Velocity averaging, kinetic formulations
 and regularizing effects in quasilinear PDEs},
Comm. Pure Appl. Math. \textbf{60} (2007), 1488--1521.

\bibitem{Ta} L. Tartar, \textit{Compensated compactness and applications to partial
differential equations}, In: Research Notes in Mathematics,
Nonlinear Analysis and Mechanics:  Herriot-Watt Symposium, Vol.
\textbf{4}, ed. R.J. Knops, Pitman Press, 1979.
\end{thebibliography}
\end{document}